\documentclass[11pt]{article}
\usepackage{amsmath}
\usepackage{xspace}
\usepackage{german}
\usepackage{amssymb}
\usepackage{graphicx}
\usepackage{epsfig}
\newtheorem{theorem}{Theorem}[section]

\newtheorem{proposition}[theorem]{Proposition}

\newtheorem{lemma}[theorem]{Lemma}

\def\IC{{\mathbf C}}

\def\0{{\mathbf{0}}}

\def\diag{{\rm diag}\,}
\def\tr{{\rm tr}\,}

\def\qed{\hfill\vbox{\hrule width 6 pt
    \hbox{\vrule height 6 pt width 6 pt}}\medskip}

\def\[{\left [}
\def\]{\right ]}
\def\({\left (}
\def\){\right )}

\def\ba{\begin{array}}
\def\barl{\begin{array}{rl}}
\def\barrl{\begin{array}{rrl}}
\def\ea{\end{array}}
\def\bt{\begin{tabular}}
\def\et{\end{tabular}}
\def\bit{\begin{itemize}}
\def\eit{\end{itemize}}
\def\beu{\begin{enumerate}}
\def\eeu{\end{enumerate}}

\newcommand{\maxover}[1]{\ensuremath{\underset{#1}\max}}
\newcommand{\minover}[1]{\ensuremath{\underset{#1}\min}}
\newcommand{\unity}{\ensuremath{{\rm 1 \negthickspace l}{}}}
\newcommand{\zero}{\ensuremath{\mathbb O}{}}

\newcommand{\braket}[2]{\ensuremath{\langle #1 | #2 \rangle}{}}

\newcommand{\grad}{\operatorname{grad}}

\renewcommand{\Re}{\operatorname {Re}}

\newcommand{\eig}{\operatorname {eig}}
\begin{document}

\title{Least-Squares Approximation by Elements from Matrix Orbits\\
Achieved by Gradient Flows on Compact Lie Groups}

\author
{Chi-Kwong Li\thanks{Department of Mathematics, College of William and Mary,
Williamsburg, VA 23185, USA. Li is an honorary professor of the University of 
Hong Kong. His research was partially supported
by USA NSF and HK RGC. e-mail:ckli@math.wm.edu}, 
Yiu-Tung Poon\thanks{Department of Mathematics, Iowa State University,
Ames, IA 50051 USA. e-mail: ytpoon@iastate.edu}, 
and Thomas Schulte-Herbr\"{u}ggen\thanks{Department of Chemistry, Technical 
University Munich (TUM), D-85747, Garching-Munich, Germany. e-mail: tosh@ch.tum.de;
T.S.H. is supported in part by the EU-programme QAP and by 
the excellence network of Bavaria through QCCC.}}

\date{dated: \today}
\maketitle

\smallskip\noindent
{\bf Abstract}

Let $S(A)$ denote the orbit of a complex or real matrix $A$ under a certain
equivalence relation such as unitary similarity, unitary equivalence, unitary congruences
etc.
Efficient gradient-flow algorithms are constructed to determine
the best approximation
of a given matrix $A_0$ by the sum of matrices
in $S(A_1), \dots, S(A_N)$ in the sense of finding the Euclidean
least-squares distance
$$\min\left\{\|X_1+ \cdots + X_N - A_0\|: X_j \in S(A_j), \ j = 1, \dots, N\right\}.$$
Connections of the results to different pure and applied areas
are discussed.



\section{Introduction}
\setcounter{equation}{0}

Motivated by problems in pure and applied areas, there has been a
great deal of interest in studying equivalence classes on matrices, say,
under compact Lie group actions. For instance,

(a) the unitary (orthogonal)
similarity orbit of a complex (real) square matrix $A$ is the set of
matrices of the form $UAU^*$ for unitary (or real orthogonal) matrices $U$,

(b) the unitary (orthogonal)
equivalence orbit of a complex (real) rectangular matrix $A$
is the set of matrices of the form $UAV$ for unitary (orthogonal)
matrices
$U, V$ of appropriate sizes,

(c) the unitary $t$-congruence orbit of a complex square matrix $A$ is
the set of matrices of the form $UAU^t$ for unitary matrices $U$,

%(c') the unitary congruence orbit under complex conjugation
%of a complex square matrix $A$ is the set of matrices of the form $UA\bar U$,
%
%{\bf I am not sure (c')  is a natural group action. CK}
%
%{\bf Actually, this is not even an equivalence relation. YT}
%
%Let
%$A= \left(\begin{smallmatrix} 1 & 0\\  0 & 0\end{smallmatrix}\right)$, $U=\left(\begin{smallmatrix} 0 & i\\ 1  & 0\end{smallmatrix}\right)$, $V=U^T$, $B=UA\overline{U}$ and $C=VB\overline{V}=UVA\overline{U}\overline{V}=\left(\begin{smallmatrix} -1 & 0\\  0 & 0\end{smallmatrix}\right)$. Suppose $W=\(w_{ij}\)$ is a unitary matrix satisfying $WA\overline{W}=C$. Then we have $WA=CW^T\Rightarrow w_{11}=w_{12}=w_{21}=0\Rightarrow W$ is not  unitary.
%
(d) the orthogonal similarity orbit of a complex square matrix
$A$ is the set of
matrices of the form $QAQ^t$ for complex orthogonal matrices $Q$,
i.e., $Q^tQ = I_n$,

(e) the similarity orbit of a square matrix $A$ is the set of matrices
of the form $SAS^{-1}$ for invertible matrices $S$.

\medskip
It is often useful to determine whether a matrix $A_0$ can be
written as a sum of matrices from orbits $S(A_1), \dots, S(A_N)$.
Equivalently, one would like to know whether
$$S(A_0) \subseteq S(A_1) + \cdots + S(A_N).$$
For $N = 1$, it reduces to the
basic problem of checking whether $A_0$ is equivalent to
$A_1$. In some cases, even this is non-trivial.
 For instance, it is not easy to check
whether two $n\times n$ complex matrices are unitarily similar.
For $N > 1$, the problem is usually more involved.
Even if there are theoretical results, it may not be easy
to use them in practice or checking examples of matrices
of moderate sizes. 
For instance,
given $10\times 10$ Hermitian matrices $A, B, C$,
to conclude that $C = UAU^* + VBV^*$ for some unitary matrices
$U$ and $V$, one needs to
check thousands of inequalities involving the eigenvalues of $A$, $B$,
and $C$; see \cite{Fu}.
Therefore, one purpose of this paper is to set up a general framework to
develop efficient computer algorithms and programs to solve such problems.
In fact, we will treat the more general problem of finding the best approximation
of a given matrix $A_0$ by the sum of matrices
from matrix orbits $S(A_1), \dots, S(A_N)$.
In other words, for given matrices $A_0, A_1, \dots, A_N$, we
determine
$$\min\left\{\|X_1 + \cdots + X_N - A_0\|:
(X_1, \dots, X_N) \in S(A_0)\times \cdots \times S(A_N)\right\}.$$
The results will be useful in solving numerical problems
efficiently, and helpful in testing conjectures of theoretical
development of the topics under considerations. As we will see in the
following discussion, some numerical examples indeed lead to
general theory; see Section 3.]

\medskip
We will consider different matrix orbits in the next few sections.
In each case, we will mention the motivation of the problems and
derive the gradient flows for the respective orbits, which will be used
to design the algorithms and computer programs to solve the optimization
problem. Note that we always consider the orbits of
similarity $SAS^{-1}$ and equivalence $SAT$, where $\left\{S,T\right\}$ can be elements
of any semisimple compact connected matrix Lie group,
in particular the special unitary group $SU(n)$
and subgroups thereof. Since these matrix Lie groups are
compact, they are themselves smooth Riemannian manifolds $M$,
which in turn implies they are endowed with a Riemannian
metric induced by the non-degenerate Killing form related to a bi-invariant
scalar product $\braket\cdot\cdot _x$ on their tangent
and cotangent spaces $T_xM$ and $T^*_xM$. The metric smoothly varies with
$x\in M$ and allows for identifying the Fr{\'e}chet differential in $T^*_xM$
with the gradient in $T_xM$. Moreover, in Riemannian manifolds the existence
and convergence of gradient flows with appropriate discretization schemes are
elaborated in detail in Ref.~\cite{SDHG07a}. In the present context, it is
important to note that the subsequent gradient flows on the unitary congruence
orbit and the unitary equivalence orbit are fundamental. The flows on compact
connected subgroups of $SU(n)$ such as $SO(n)$ or $SU(2)^{\otimes m}$ (with
$2^m=n$) can readily be derived from the flows on $SU(n)$
\cite{SDHG07,SDHG07a}.
Furthermore, in each case, we will provide numerical examples
to illustrate their efficiency and accuracy.

The situation in the general linear group $GL(N)$ and its subgroups
that are not in the intersection with the unitary groups is entirely
different: those groups are no longer
compact, but only locally compact.
For $GL(N)$~orbits we give an outlook with
some analytical results in infinma of Euclidean distances.
Since locally compact Lie groups lack {\em bi-invariant}
metrics on the tangent spaces to
their orbit manifolds, they can only be endowed with left-invariant {\em or}
right-invariant metrics. Moreover, the exponential map onto locally
compact Lie groups is no longer geodesic as in the compact case.
Consequently, one will have to devise other approximations
to the respective geodesics than obtained by the (Riemannian)
exponential. These numerics
are thus a separate topic of current research and will therefore
be pursued in a follow-up study.

With regard to notation, unless stated otherwise, the norm $||A||$
shall always be
read as Frobenius norm $||A||_2 := \sqrt{\tr\left\{A^* A\right\}}$.

\section{Unitary Similarity Orbits}
\subsection{The Hermitian Matrix Case}

For an $n\times n$ Hermitian matrix $A$, let
$S(A)$ be the set of matrices unitarily similar to $A$.
Then
$$S(A)+S(B) = \left\{X+Y: (X,Y) \in S(A)\times S(B)\right\}$$
is a union of unitary similarity orbits. Researchers have
determined the necessary and sufficient conditions of $S(C)$ to
be a subset of $S(A)+S(B)$ in terms of the eigenvalues of
$A, B$ and $C$; \cite{CW1,CW2,DST,Fu,H2,K,TF1,TF2}.
In particular, suppose $A, B, C$ have eigenvalues
$$a_1 \ge \cdots \ge a_n, \qquad
b_1 \ge \cdots \ge b_n, \qquad \hbox{ and } \qquad
c_1 \ge \cdots \ge c_n,$$
respectively. Then $S(C) \subseteq S(A) + S(B)$ if and only if
\begin{equation}\label{tr}
\sum_{j=1}^n (a_j+b_j - c_j) = 0\end{equation} and a collection of
inequalities in the form
\begin{equation}\label{rst}\sum_{r \in R} a_r +
\sum_{s\in S} b_s \ge \sum_{t \in T} c_t\end{equation} for certain
$m$ element subsets $R,S,T \subseteq \left\{1, \dots, n\right\}$  with
$1 \le m < n$ determined by the Littlewood-Richardson rules;
see \cite{DST,Fu} for details.
The study
has connections to many different areas such as representation
theory, algebraic geometry, and algebraic combinatorics, etc.
Note that the relation between Horn's problem and the Littlewood-Richardson rules
has recently also attracted attention in quantum information \cite{C08}.
The set of inequalities in (\ref{rst}) grows exponentially with $n$.
Therefore, it is not easy to check the conditions even for a moderate
size problem, say, for  $10\times 10$ Hermitian matrices.
As a matter of fact, the theory has been extended to determine
whether $S(A_0)$ is a subset of $S(A_1) + \cdots + S(A_N)$
for given  $n\times n$ Hermitian matrices $A_0, \dots, A_N$,
in terms of equality and linear inequalities
of the eigenvalues of the given matrices.
Of course, the number of inequalities involved are more
numerous.  There does not seem to be an efficient way
to use these results in practise or testing numerical examples
or conjecture in research.

It is interesting to note that by the saturation conjecture (theorem)
(see \cite{BF} and its references),
there exist Hermitian matrices with
nonnegative integral eigenvalues $a_1 \ge \cdots \ge a_n$,
and $b_1 \ge \cdots \ge b_n$
such that $A+B$ has nonnegative integral eigenvalues
$c_1 \ge \cdots \ge c_n$ if and only if
the Young diagram corresponding to
$(c_1, \dots, c_n)$ can be obtained from those of
$(a_1, \dots, a_n)$ and $(b_1, \dots, b_n)$.

\subsection{The General Complex Matrix Case}

Likewise, we study the problem
$$\min\left\{\| \sum_{j=1}^N U_j A_j U^*_j - A_0 \|:
U_1, \dots, U_N \in SU(n)\; \hbox{unitary} \right\}$$
for general complex matrices $A_0, \cdots A_N$.
Even for $N = 1$, the result is highly nontrivial.
In theory, it is related to the problem of determining
whether $A_0$ and $A_1$ are unitarily similar; see \cite{Sh}.
Also, to determine
$$\min\left\{\|UAU^* - C^* \|:
U \hbox{ unitary} \right\}$$
for $A, C \in M_n$
leads to the study of the $C$-{\it numerical range} and the
$C$-{\it numerical radius} of $A$ defined by
$$W(C,A) = \left\{ \tr(CUAU^*): U \in SU(n) \right\},$$
and
$$r(C,A) = \max\left\{|\mu|: \mu \in W(C,a)\right\}.$$
The $C$-numerical radius is important in the study of
{\it unitary similarity
invariant} norms on $M_n$, i.e., norms $\nu$ satisfy
$\nu(UXU^*) = \nu(X)$ for all $X, U \in M_n$ such that $U$ is unitary.
For instance, it is known that for every unitary similarity invariant
norm $\nu$ there is a compact subset $S$ of $M_n$ such that
$$\nu(X) = \max\left\{ r(C,X): C \in S\right\}.$$
So, the $C$-numerical radii can be viewed as the building blocks
of unitary similarity invariant norms. We refer readers to
the survey \cite{Li} for further results on the $C$-numerical range and
$C$-numerical radius. For applications of
$C$-numerical ranges in quantum dynamics,
see also Ref.~\cite{SDHG07}

For two matrices, one may study whether
$C = UAU^* + VBV^*$
for, e.g., a Hermitian $A$ and a skew-Hermitian $B$. In other words,
we want to study whether a matrix can be written as the
sum of a Hermitian matrix and a skew-Hermitian matrix
with prescribed eigenvalues.

%%%%%%%%%%%%%%%
\subsection{Sum of Hermitian and Skew-Hermitian Matrices}
%%%%%%%%%%%%%%%%%%%%
%\begin{figure}[Ht!]
%\begin{center}
%\includegraphics[width=0.65\textwidth]{LiHeSk.eps}
%\end{center}
%\caption{\label{fig:HeSk}
%Experiment 1 shows that
%(a) $\min\left\{ \|UAV^* + XBY^* -  C\|\;:\; U,V,X,Y \in SU(n)\right\} = 0$ while
%(b) $\min\left\{ \|U'AU'^* + V'BV'^* -  C\|\;:\; U',V' \in SU(n)\right\} > 0$
%for matrix dimensions $n=2,5,10,20,50$.
%}
%\end{figure}
%%%%%%%%%%%%%%%%%%%%

For $C = UAU^* + VBV^*$ with $A = A^*$ and $B = -B^*$, there are
many known inequalities relating the eigenvalues of $A$ and $B$
to the eigenvalues and singular values of $C$; see \cite{CHL} and
the references therein. However, there has been no known necessary and
sufficient
condition for the existence of matrices $A, B, C$ satisfying
$C = UAU^*+VBV^*$ with $A = A^*$ and $B=-B^*$ with prescribed
eigenvalues or with prescribed singular values.
Nevertheless, it is easy to solve the approximation
problem
$$\min\left\{\|U^*AU + V^*BV - C\|: U, V \hbox{ unitary} \right\}.$$
The following result actually holds for any {\it unitarily invariant}
norm on $n\times n$ matrices using the same proof; see \cite{LT}.
Furthermore, we can use this result to verify that our algorithm indeed
yield the optimal solution; see Example 2 in Section 2.5.

%\medskip
%\noindent
%{\bf Experiment 1:}\\
%We generate a (random) upper triangular $\tilde C \in M_n$ and set
%$A := (\tilde C+\tilde C^*)/2$ and $B := (\tilde C-\tilde C^*)/(2i)$.
%Next, we change the strictly upper triangular part of $\tilde C$ to get
%$C$ so that $C$ and $\tilde C$ have the same eigenvalues.
%Fig.~\ref{fig:HeSk}(a) shows that $C$ is on the sum of
%the two unitary equivalence orbits of $A$ and $B$ respectively,
%whereas Fig.~\ref{fig:HeSk}(b) demonstrates it is {\em not}
%on the sum of the two corresponding similarity orbits.
%
%{\small
%Check whether $C = XAY + RBS$ using the results/programs
%in Section 3. If yes, then all the singular inequalities
%relating $A,B$ and $C$ hold.
%Since $C$ and $\tilde C$ have the same eigenvalues, all possible
%eigenvalue inequalities relating those of $C$ and those of
%$A, B$ should also hold.
%
%Now, if we can show that $\min\| C - UAU^* - VBV^*\| > 0$,
%it will illustrate that the eigenvalues and singular valaues
%inequalities are not enough to ensure the existence of
%$C = UAU^* + VBV^*$ with $A = A^*$ and $B = -B^*$.
%{\bf YES, THIS IS THE CASE}
%
%On the other hand, if many examples support that
%$\min\|C - UAU^* - VBV^*\| = 0$ once the singular values
%inequalities hold (and the eigenvalues inequalities hold automatically
%by our construction), then we have more confidence that
%combining these inequalites will give the necessary
%and sufficient condition.
%{\bf NO, NOT A SINGLE INSTANCE FOUND TO SUPPORT THIS}
%}

\begin{theorem}   Let $\|\cdot\|$ be the Frobenius norm on $M_n$.
Let $A, B, C \in M_n$ with $A = A^*$ and $B = -B^*$.
Suppose $U, V \in M_n$ are unitary matrices
such that $U\tfrac{1}{2}(C+C^*)U^* = \diag(f_1, \dots, f_n)$ with
$f_1 \ge \cdots \ge f_n$,
and $V\tfrac{1}{2}(C-C^*)V^* = i\; \diag(g_1, \dots, g_n)$
with $g_1 \ge \cdots \ge g_n$.
Suppose $A$ is unitarily similar to a diagonal matrix
$A_1$ (respectively, $A_2)$ with  diagonal entries arranged in
descending (respecitively, ascending) order.
Suppose $-i B$ is unitarily similar to a diagonal matrix
$-i B_1$ (respectively, $-i B_2)$ with  diagonal entries arranged in
descending (respecitively, ascending) order.
Then
\begin{eqnarray*}
\|U^*A_1U + V^*B_1V -C \|^2 &=& {\sum_{j=1}^n(|f_j - a_j|^2 + |g_j-b_j|^2)}\\
\|U^*A_2U + V^*B_2V -C\|^2 &=& {\sum_{j=1}^n(|f_j - a_{n-j+1}|^2 + |g_j-b_{n-j+1}|^2)}
\end{eqnarray*}
and for any unitary $X, Y \in M_n$,
$$\|U^*A_1U + V^*B_1V -C\| \le \|X^*AX +  Y^*BY -C\|
\le \|U^*A_2U + V^*B_2V -C\|.$$
\end{theorem}

\it Proof. \rm
Let $F = \tfrac{1}{2}(C+C^*)$ and $G = \tfrac{-i}{2}(C-C^*)$.
It is well known that
$$\|F - U^*A_1U\| \le \|F - X^*AX\| \le \|F - U^*A_2U\|$$
and
$$\|G - V^*B_1V\| \le \|G - Y^*BY\| \le \|G - V^*B_2V\|$$
for any unitary $X, Y \in M_n$; see \cite{LT}.
Since $\|H+iK\|^2 = \|H\|^2 + \|K\|^2$
for any Hermitian $H,K \in M_n$, the results follow.
\qed

%%%%%%%%%%%%%%%

%%%%%%%%%%%%%%%
\subsection{Deriving Gradient Flows on Unitary Similarity Orbits}\label{sec:flow_similarity}
%%%%%%%%%%%%%%%

To begin with, we focus on the problem of approximating a given
matrix $C$ using matrices from two unitary similarity orbits,
i.e., finding
$$\min \left\{ \|UAU^* + VBV^* - C\|: U, V \in SU(n)\; \hbox{ unitary} \right\}.$$
For simplicity, here we describe the steepest descent method to
search for unitary matrices $U_0, V_0$ attaining the optimum.
Refined approaches like conjugate gradients, Jacobi-type or Newton-type methods
may be implemented likewise, see for instance \cite{SDHG07a}.
As will be shown below, more than two unitary similarity orbits can be treated
similarly.
The basic idea is to improve the current unitary pair $(U_k,V_k)$ to
$(U_{k+1}, V_{k+1})$ so that
$$\|U_{k+1}A U^*_{k+1} + V_{k+1}B V^*_{k+1} - C\| <
\|U_{k}A U^*_{k} + V_{k}BV^*_{k} - C\|$$
until the successive iterations differ only by a small tolerance,
or the gradient ({\em vide infra}) vanishes.
Further, to avoid pitfalls by local minima whenever the Euclidean distance
cannot be made zero, we use a sufficiently large multitude
of different random starting points $(U_0,V_0)$ for our algorithm.
Needless to say, a positive matching result is constructive,
while a negative result may be due to local minima. It is therefore
important to use a sufficiently large set of initial conditions
for confident conclusions in the negative case.

For a start, consider the least-squares minimization task
\begin{equation}
\minover{U,V\in SU(n)} ||UAU^* + VBV^* - C||_2^2\;,
\end{equation}
which can be rewritten as
\begin{eqnarray*}
%\begin{split}
||UAU^* &+& VBV^* - C||_2^2 \\
&=& ||UAU^* + VBV^*||_2^2 + ||C||_2^2 - 2 \Re \tr\left\{C^*(UAU^* + VBV^*)\right\} \\
&=& ||A||_2^2 + ||B||_2^2 + ||C||_2^2 -
2 \Re \tr\left\{C^*(UAU^* + VBV^*) - UAU^*\; VB^*V^*\right\}
%\end{split}
\end{eqnarray*}
and thus is equivalent to the maximisation task
\begin{equation}
\maxover{U,V\in SU(n)} \Re \tr\left\{C^*(UAU^* + VBV^*) - UAU^*\; VB^*V^*\right\}\;.
\end{equation}

\noindent
Therefore we set
\begin{equation}
f(U,V) := \tr\left\{(UAU^* + VBV^*)\, C^* - UAU^*\; VB^*V^*\right\}
\end{equation}
and $F(U,V) := \Re f(U,V)$.
Then its Fr{\'e}chet derivative $D_U f(U): T_U\mathcal U \to T_{f(U)}\mathcal
U$ can be seen as a tangent map, where the elements of the tangent space
$T_U\mathcal U$ to the Lie group of unitaries $\mathcal U=SU(n)$ or $U(n)$ at the point $U$
take the form $\Omega U$ with $\Omega=-\Omega^*$ being itself an element of
the Lie algebra. The differential thus reads
\begin{eqnarray*}
%\begin{split}
D_U f(U)(\Omega U)
&=& \tr\left\{((\Omega U) A U^* + UA(\Omega U)^*)(C^* - VB^*V^*)\right\} \\
&=& \tr\left\{((\Omega U) A U^* - UAU^*(\Omega U) U^*)(C^* - VB^*V^*)\right\}\\
&=& \tr\left\{( A U^* (C^* - VB^*V^*) - U^*(C^* - VB^*V^*)\, UAU^*) (\Omega U)\right\}
%\end{split}
\end{eqnarray*}
where we used the invariance of the trace under cyclic permutations and
$(\Omega U)^* = - U^* (\Omega U) U^*$, which follows from the product rule for
$D(\unity)(\Omega U) = D(UU^*)(\Omega U) =0 = (\Omega U)U^* + U(\Omega U)^*$
in consistency with the Lie-algebra elements $\Omega$ being skew-Hermitian.
Moreover, by identifying
\begin{equation}
D_U f(U)\cdot (\Omega U)
= \braket{\grad_U f(U)}{\Omega U} = \tr\left\{(\grad_U f(U))^* \Omega U\right\}
\end{equation}
one finds
\begin{equation*}
\grad_U f(U)  = (C - VBV^*) U A^* - UA^*U^* (C - VBV^*) U
           = \big[(C - VBV^*), UA^*U^*\big]\,U\quad.
\end{equation*}
With $[X^*,Y]_s
:= \tfrac{1}{2}([X^*,Y] - [X^*,Y]^*) = \tfrac{1}{2}([X^*,Y] + [X,Y^*])$
as skew-hermitian part of the commutator
one obtains for $F(U):= \Re f(U)$
\begin{equation}
\grad_U F(U) = \big[(C^* - VB^*V^*), UAU^*\big]_s\,U\quad.
\end{equation}
Taking the respective Riemannian exponentials
$\exp_U(\grad_U F(U))$ and $\exp_V(\grad_V F(V))$ thus gives
the recursive gradient flows
\begin{eqnarray*}
%\begin{split}
U_{k+1} &=& \exp\left\{-\alpha_k [ U_kAU_k^*, (C^* - V_kB^*V_k^*)]_s\right\}\;U_k\\
V_{k+1} &=& \exp\left\{-\beta_k [ V_kBV_k^*, (C^* - U_kA^*U_k^*)]_s\right\}\;V_k
%\end{split}
\end{eqnarray*}
as discretized solutions of the coupled gradient system
\begin{equation}
\dot U = \grad_U F(U,V) \quad\text{\rm and}\quad \dot V = \grad_V F(U,V) \;.
\end{equation}
Conditions for convergence are described in detail in \cite{HM94}.
For appropriate step sizes $\alpha_k,\beta_k$ see also Ref.~\cite{NMRJOGO}.

Generalizing the findings from a sum of two orbits to
higher sums of unitary orbits is straightforward: the problem
\begin{equation}
\min\left\{\| \sum_{j=1}^N U_j A_j U^*_j - A_0 \|:
U_1, \dots, U_N \in SU(n)\;\hbox{unitary} \right\}
\end{equation}
can be addressed by the system of coupled gradient flows ($j=1,2,\dots, N$)
\begin{equation}\label{eqn:US-flows}
U_{k+1}^{(j)} = \exp\left\{-\alpha_k^{(j)}
{[ A_k^{(j)}, A_{0jk}^*]}_s\right\}\;U_k^{(j)}
\end{equation}
where for short we set $A_k^{(j)} := U_k^{(j)}A_j{U_k^{(j)}}^*$ and
$A_{0jk}:= A_0 - \sum\limits_{\nu=1 \atop \nu\neq j}^N A_k^{(\nu)}$.

These gradient flows follow the extension of the original idea on the
orthogonal group \cite{Bro88,HM94} to the unitary group \cite{Sci98}, where
here we introduce a larger system of coupled flows.

\subsection{Numerical Examples}

Here we demonstrate gradient flows minimising
$\| \sum_{j=1}^N U_j A_j U_j^* - A_0 \|$ over the unitaries $ U_1, \dots, U_N$
for given Hermitian matrices $A_0, \cdots A_N$.

\medskip
\noindent
{\bf Example 1}\\
As a test case, consider the following examples for
finding $U_j \in \IC^{10\times 10}$. For $j=1,2,\dots,N$
choose a set of random unitaries $U_j^{(r)}\in\IC^{10\times 10}$
distributed according to the Haar measure as recently described in \cite{Mez07}
and define
$A_j := \diag(1,3,5,\dots,19) +
\tfrac{j-1}{10}\unity_{10}$ and $
	A_0^{(N)}:= \diag(a_1, ..., a_{10})$
where $a_1,a_2,\dots,a_{10}$ are the eigenvalues of
$A'_{0,N}:= \sum_{j=1}^N U_j^{(r)}
A_j {U_j^{(r)}}^*$ (and $\unity_{10}$ is the $10\times 10$ unity matrix).
%%%%%%%%%%%%%%%%%%%
\begin{figure}[Ht!]
\begin{center}
\includegraphics[width=0.65\textwidth]{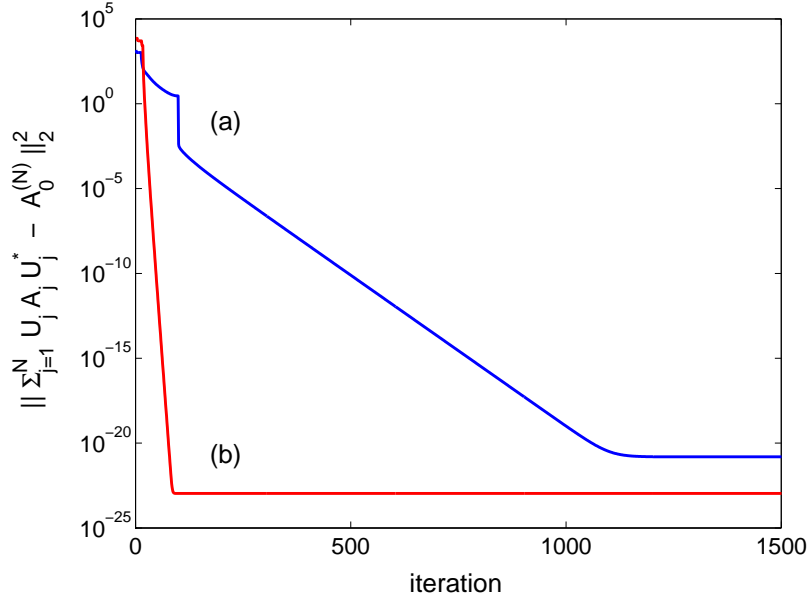}
\end{center}
\caption{\label{fig:UAU_N}
Coupled flows minimizing $||\sum_{j=1}^N U_jA_jU_j^* - A_0^{(N)}||^2_2$ with
(a) $N=2$ and (b) $N=10$ for Example 1.
}
\end{figure}
%%%%%%%%%%%%%%%%%%%
As shown in Fig.~\ref{fig:UAU_N}, the gradient flow of Eqn.~\ref{eqn:US-flows} minimizes
$||\sum_{j=1}^N U_jA_jU_j^* - A_0^{(N)}||^2_2$ by driving it practically to
zero. Note that in Fig.~\ref{fig:UAU_N}b the combined flow on $N=10$ unitaries
converges even faster than in Fig.~\ref{fig:UAU_N}a, where $N=2$ and the flow
is more sensitive to saddle points as may be inferred from the jumps in trace
(a).

%%%%%
%%%%%
\bigskip
\noindent
{\bf Example 2}\\
Let $A,B$ be Hermitian and $C$ arbitrary, e.g.,
$ A = \left(\begin{smallmatrix} 2 & 5 & 11\\  5 & 8  & 15\\  11 & 15 & 16\end{smallmatrix}\right), \;
 B = \left(\begin{smallmatrix} 6 & 8 & 9\\  8 & 12 & 10\\  9  & 10 & 0\end{smallmatrix}\right), \;
 C = \left(\begin{smallmatrix} 1 & 11 & 3\\ 6 & 9 & 3\\ 8 & 9 & 2\end{smallmatrix}\right) \;.
$
Then $a := \eig(A)= (-5.6674; -0.4830; 32.1504),
b := \eig(B)= (-7.4816; 0.7123; 24.7693)$ and
$f: = \eig\tfrac{1}{2}(C+C^*)=(-4.9555; -1.3888; 18.3443), g := \eig\tfrac{-i}{2}(C-C^*)=(-4.6368; 0; 4.6368)$.
According to Theorem 2.1 one gets
\begin{equation}
\Delta:= \minover{U,V\in SU(3)} ||UAU^*+iVBV^*-C||_2^2 = (a-f)^*(a-f) + (b-g)^*(b-g) = 605.8521\;.
\end{equation}
More precisely $\Delta = 605.852131091'3004$, while 100 runs of the gradient flow with independent
random initial conditions give a mean $\pm$ rmsd.~of $\bar \Delta = 605.852131091'3570 \pm 1.13\cdot 10^{-10}$.
%%%%%
%%%%%

\section{Unitary Equivalence}

In this section, we study
$$\min\left\{\|\sum_{j=1}^N U_jA_jV_j - A_0\|: U_1, \dots, U_N \in U(n)
{\quad \text {and}\quad } V_1, \dots, V_N \in U(m)\; \hbox{ unitary} 
\right\}$$
for rectangular matrices $A_0, \dots, A_N$.
By the result of O'Shea and Sjamaar \cite{OS},
$$\min\|\sum_{j=1}^N U_jA_jV_j - A_0\|=0$$
if and only if
$$\min\|\sum_{j=1}^N W_j^*\tilde A_jW_j - \tilde A_0\|=0$$
where
\begin{equation*}
\tilde A_j = \begin{pmatrix} 0 & A_j \cr
A_j^* & 0 \cr\end{pmatrix} \qquad \hbox{ for } j=0,1, \dots, N.
\end{equation*}
Thus, by the results concerning unitary similarity orbits (see Section 2),
\begin{equation}\label{eq1}
\min\left\{\|A_0 - \sum_{j=1}^N U_jA_jV_j\|:
U_1, \dots, U_N;\, V_1, \dots, V_N \; \hbox{ unitary} \right\} = 0
\end{equation}
if and only if the singular values of $A_0, A_1, \dots, A_N$
satisfy a certain set of linear inequalities.
Clearly,
$\min\left\{\|A - UBV\|: U, V \hbox{ unitary} \right\} = 0$
if and only if $A$ and $B$ have the same singular values.
In general,
it is interesting to check whether
$$\sqrt 2 \min\|\sum_{j=1}^N U_jA_jV_j - A_0\|
= \min\|\sum_{j=1}^N W_j^*\tilde A_jW_j - \tilde A_0\|=0.$$
In computer experiments
(see Example 6 in Section 3), we observe that
(\ref{eq1}) always holds if $A_0, A_1, \dots, A_N$ are randomly
generated matrices generated by {\sc matlab}.
We explain this phenomenon in the following.
We begin with a simple observation.

\begin{lemma}
Suppose $a_0, a_1, \dots, a_N \in (0, \infty)$.
The following are equivalent.

{\rm (a)} There are complex units $e^{it_1}, \dots, e^{it_N}$ such that
$a_0 - \sum_{j=1}^N  a_j e^{it_j} = 0.$

{\rm (b)} There is an $N+1$ side convex polygon whose sides have lengths
$a_0, \dots, a_N$.

{\rm (c)} $\sum_{j=0}^N a_j - 2a_k \ge 0$ for all $k = 0, 1, \dots, N$.

\end{lemma}

Form this observation, one easily gets the following
condition related to the equality (\ref{eq1}).

\begin{proposition} \label{prop1} Let $A_j = \diag(a_{1j}, \dots, a_{nj})$ be
nonnegative diagonal matrices for $j = 0, 1, \dots, N$, and let
$v_j = (a_{1j}, \dots, a_{nj})^t$.
Then there exist permutation matrices $P_1, \dots, P_N$
and diagonal unitary matrices $D_1, \dots, D_N$ such that
$$A_0 = \sum_{j=1}^N D_jP_jA_jP_j^t$$
if and only if the entries of each row of the matrix
$$[v_0 | P_1v_1 | \cdots | P_N v_N]$$
correspond to the sides of a $N+1$ side convex polygon.
\end{proposition}

\noindent
If one examines the singular values of an $n\times n$
random matrix generated by {\sc matlab}, we see that there is
always a dominant singular values of size about $n/2$, and
the other singular values range from 0 to $1.5n$ in a rather
systematic pattern.  So, it is often
possible to apply Proposition \ref{prop1} to get equality
(\ref{eq1}) if $A_0, \dots, A_N$
are random matrices generated by {\sc matlab} for $N \ge 2$.

\medskip
In contrast, for general matrices, it is easy to construct
$A_0, A_1, \dots, A_N$ such that (\ref{eq1}) fails.

\medskip
\medskip
\noindent
{\bf Example 3}\\
Let $A_0 = \diag(N^2,N+1) \oplus 0_{n-2}$
and $A_j = \diag(N,1) \oplus 0_{n-2}$ for $j = 1, \dots, N$.
Then clearly Eqn.~\ref{eq1} does not apply, because
$$\sum_{j=1}^n s_j(A_0) > \sum_{i=1}^N \sum_{j=1}^n s_j(A_j).$$

%%%%%%%%
%%%%%%%%
%%%%%%%%%%%%%%%%%%%
%\begin{figure}[Ht!]
%\begin{center}
%\includegraphics[width=0.65\textwidth]{Ex3-3.eps}
%\end{center}
%\caption{\label{fig:Ex3-3}
%Illustrative verification of Example 3.3 --- internal for the authors.
%Here $N=3$ and $n-2 =4$.
%}
%\end{figure}
%%%%%%%%%%%%%%%%%%%%
%%%%%%%%
%%%%%%%%
Recall that the Ky Fan $k$-norm of a matrix $A \in M_n$
is defined as $\|A\|_k = \sum_{j=1}^k s_j(A)$,
and a norm $\|\cdot\|$ on $M_n$ is unitarily invariant if
$\|A\| = \|UAV\|$ for all $A \in M_n$ and unitary $U, V \in M_n$.
By the Ky Fan dominance theorem, two matrices
$A, B \in M_n$ satisfy  $\|A\|_k \le \|B\|_k$ for $k = 1, \dots, n$
if and only if $\|A\| \le \|B\|$ for all unitarily invariant
norms $\|\cdot\|$.
In view of this example, we have the following result.

\begin{proposition} Suppose $A_0, A_1, \dots, A_N \in M_n$ satisfy
(\ref{eq1}). Then for all unitarily  invariant norms,
$$2\|A_i\| \le \sum_{j=0}^N \|A_j\|, \qquad i = 0, 1, \dots, N,$$
and equivalently, for $k = 1, \dots, n$,
\begin{equation} \label{eq2}
2\|A_i\|_k \le \sum_{j=0}^N \|A_j\|_k, \qquad i = 0,1, \dots, N.
\end{equation}
Moreover, if there is $k$ such that equality (\ref{eq2}) holds, then
(\ref{eq1}) holds if and only if
$A_j$ is unitarily similar to
$B_j \oplus C_j$ with $B_j \in M_k$ for $j=0, \dots, N$ such that
$$\min\left\{\|B_0 - \sum_{j=1}^N U_jB_jV_j\|:
U_1, \dots, U_N, V_1, \dots, V_N \in M_k \hbox{ are unitary} \right\} = 0$$
and
$$\min\left\{\|C_0 - \sum_{j=1}^N X_jC_jY_j\|:
X_1, \dots, X_N, Y_1, \dots, Y_N \in M_{n-k}
\hbox{ are unitary} \right\} = 0.$$
\end{proposition}

It would be nice if one can get (\ref{eq1}) by checking the relatively
easy condition (\ref{eq2}).
Unfortunately, the following example shows that it is not true.

\medskip
\medskip
\noindent
{\bf Example 4}\\
Let $A_0 = \diag(14,2)$, $A_1 = \diag(8,0)$, $A_2 = \diag(7,4)$.
Then (\ref{eq2}) is satisfied for all $k\ge 1$ but by the result in \cite{LP},
$$\diag\(U_1A_1V_1+U_2A_2V_2\)\ne  (14,2)$$
for all unitaries $U_i,\ V_j$.

%%%%%%%%
%%%%%%%%
%%%%%%%%%%%%%%%%%%%%
%\begin{figure}[Ht!]
%\begin{center}
%\includegraphics[width=0.65\textwidth]{Ex3-5.eps}
%\end{center}
%\caption{\label{fig:Ex3-5}
%Illustrative verification of Example 3.5 --- internal for the authors.
%}
%\end{figure}
%%%%%%%%%%%%%%%%%%%
%%%%%%%%
%%%%%%%%%%%%
\subsection{Deriving Gradient Flows on Unitary Equivalence Orbits}\label{sec:flow_equivalence}

For minimizing $||UAV - C||_2^2$ one has to maximize
$$F(U,V) := \Re \tr \left\{UAVC^*\right\} = \tfrac{1}{2} \tr \left\{UAVC^* + (UAVC^*)^*\right\}\;.$$
By the same arguments as before, from its Fr{\'e}chet differential
$$ D_U F(U,V)(\Omega U) =
\tfrac{1}{2} \tr \left\{(\Omega U) AVC^* - CV^*A^*U^* (\Omega U) U^*\right\}
= \tfrac{1}{2} \tr \left\{(AVC^* - U^*CV^*A^*U^*) (\Omega U)\right\}$$
one obtains the gradient---where henceforth we keep writing
$(\cdot)_s$ for the skew-Hermitian part
$$ \grad_U F(U,V) = \tfrac{1}{2} (AVC^* - U^*CV^*A^*U^*)^*
= -(UAVC^*)_s \;U\;. $$
An analogous result follows for $\grad_V F(U,V)$.
Taking again the respective Riemannian exponentials leads
to the recursive scheme
\begin{eqnarray*}
U_{k+1} &=& \exp\left\{-\alpha_k (U_kAV_kC^*)_s\right\}\;U_k\\
V_{k+1} &=& \exp\left\{-\beta_k (V_kC^*U_kA)_s\right\}\;V_k\;,
\end{eqnarray*}
which also can be used, {\em e.g.}, for a singular-value decomposition of $A$
by choosing $C$ real diagonal.

Likewise, minimizing $||UAV + XBY - C||^2_2$ by maximizing
$\Re \tr \left\{UAV(C-XBY)^* + XBY C^*\right\}$ translates into the same flows when
substituting $C\mapsto (C-X_kBY_k)$
%\begin{eqnarray*}
%U_{k+1} &=& \exp\left\{-\alpha_k (U_kAV_k(C-X_kBY_k)^*)_s\right\}\;U_k\\
%V_{k+1} &=& \exp\left\{-\beta_k (V_k(C-X_kBY_k)^*U_kA)_s\right\}\;V_k
%\end{eqnarray*}
with analogous recursions for $X_{k+1}$ and $Y_{k+1}$.
Along these lines, it is straightforward to address the general task
\begin{equation}
\min\left\{\|\sum_{j=1}^N U_jA_jV_j - A_0\|: U_1, \dots, U_N \in U(n)
{\quad \text {and}\quad } V_1, \dots, V_N \in U(m)\; \hbox{ unitary} 
\right\}
\end{equation}
with rectangular matrices $A_0, \dots, A_N$ by a system of
$2 N$ coupled gradient flows ($j=1,2,\dots, N$)
\begin{eqnarray}\label{eqn:UE-flows}
U_{k+1}^{(j)}
&=& \exp\left\{-\alpha_k^{(j)} (U_k^{(j)} A_j V_k^{(j)} A_{0jk}^*)_s\right\}\;U_k^{(j)}\\
V_{k+1}^{(j)}
&=& \exp\left\{-\beta_k^{(j)} (V_k^{(j)} A_{0jk}^* U_k^{(j)}A_j)_s\right\}\;V_k^{(j)}
\end{eqnarray}
where we use the short-hand
$A_{0jk}:= A_0 - \sum\limits_{\nu=1\atop \nu\neq j}^N U_k^{(\nu)} A_\nu V_k^{(\nu)}$.
\vskip .5in
%%%%%%%%%%

\subsection{Numerical Examples}
Using the flows derived in section~\ref{sec:flow_equivalence},
in this section, we study
$$\min\left\{\|\sum_{j=1}^N U_jA_jV_j - A_0\|: U_1, \dots, U_N \in U(n)
{\quad \text {and}\quad } V_1, \dots, V_N \in U(m) \hbox{ unitary} \right\}$$
for rectangular matrices $A_0, \dots, A_N$.

\medskip
\noindent
{\bf Example 5}\\
As an example of rectangular $A_j \in \IC^{10\times 15}$, consider the
analogous flows. In order to obtain $U_j \in \IC^{10\times 10}$ and $V_j \in
\IC^{15\times 15}$ for $j=1,2,\dots,N$ choose a set of random unitary pairs
$(U_j^{(r)},V_j^{(r)})\in\IC^{10\times 10} \times \IC^{15\times 15}$ and
define
$$A_j := [\, \diag(1,3,5,\dots,19)
+ \tfrac{j-1}{10}\unity_{10}\; |\; \zero_{10,5}\,] \quad\text{and}\quad%
        A_0^{(N)}:= [\, \diag(s_1, ..., s_{10})\; |\; \zero_{10,5}\,] $$
where $s_1,s_2,\dots,s_{10}$ are now the singular values of
$A'_{0,N}:= \sum_{j=1}^N U_j^{(r)} A_j V_j^{(r)}$
and $\zero_{10,5}$ is the $10\times 5$ zero-matrix.
%%%%%%%%%%%%%%%%%%%
\begin{figure}[Ht!]
\begin{center}
\includegraphics[width=0.65\textwidth]{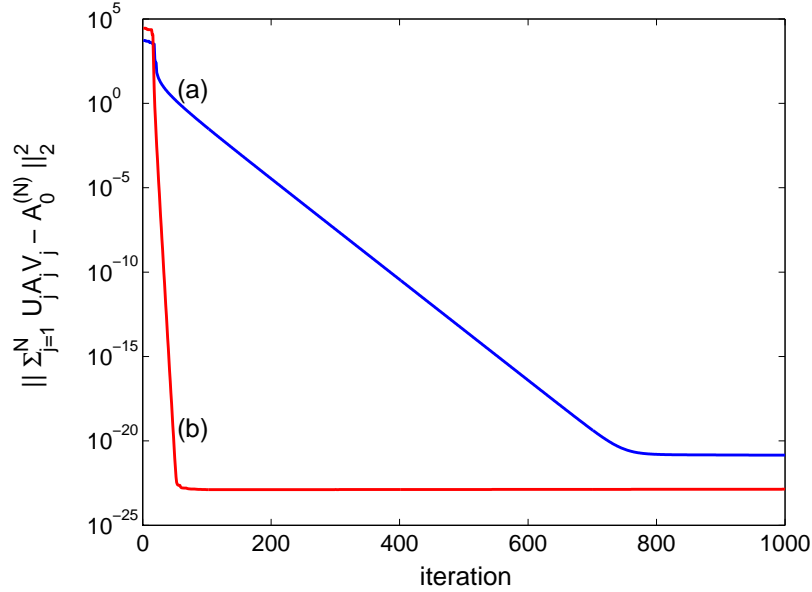}
\end{center}
\caption{\label{fig:UAV_N}
Coupled flows minimizing $||\sum_{j=1}^N U_jA_jV_j - A_0^{(N)}||^2_2$ with
(a) $N=2$ and (b) $N=10$.
Here the $A_j \in \IC^{10\times 15}$ are rectangular so that
$U_j \in \IC^{10\times 10}$ and $V_j \in \IC^{15\times 15}$.
}
\end{figure}
%%%%%%%%%%%%%%%%%%%
Fig.~\ref{fig:UAV_N} shows how the coupled gradient flow
minimizes $||\sum_{j=1}^N U_jA_jV_j - A_0^{(N)}||^2_2$
by driving it practically to zero.
Again the combined flow on $N=10$ unitary pairs
(Fig.~\ref{fig:UAV_N}b) converges faster
than the one for $N=2$ unitary pairs given in Fig.~\ref{fig:UAV_N}a.

\subsubsection{Observation Concerning Sums of Unitary Equivalence Orbits}
\begin{figure}[Ht!]
\begin{center}
\includegraphics[width=0.65\textwidth]{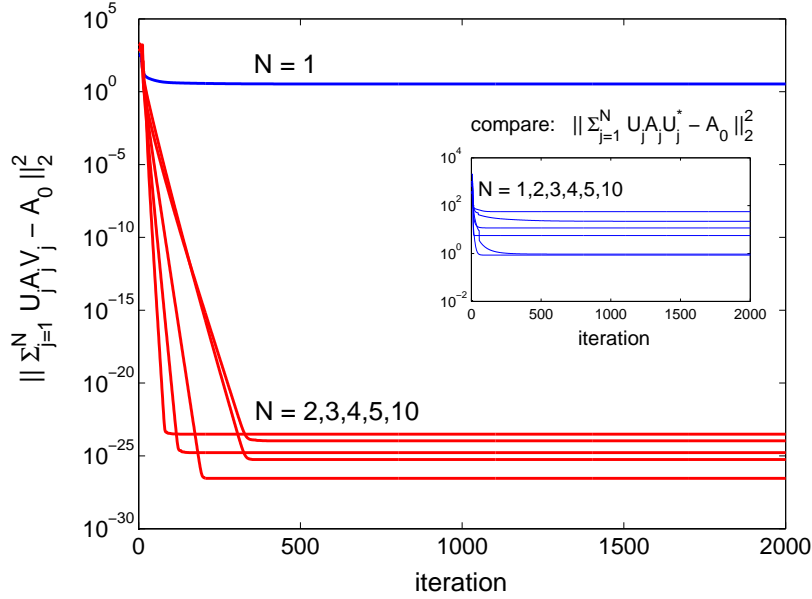}
\end{center}
\caption{\label{fig:UAV_N2}
A random complex square matrix $A_0\in \IC^{10\times 10}$ is typically distant
from a single ($N=1$) equivalence orbit
of another random square matrix $U A_1 V$, as shown in the upper trace.
However, it is typically arbitrarily close to a sum of equivalence
orbits of several independent random square matrices as demonstrated
in the lower traces:
$\|\sum_{j=1}^N U_jA_jV_j - A_0\|^2_2 \to 0$ for $N=2,3,4,5,10$.
In contrast, the inset shows this does not hold for $N=1$ through $N=10$
for similarity orbits $\sum_{j=1}^N U_jA_jU^*_j$.
}
\end{figure}

A non-zero random complex matrix $A_0$ is typically distant from a single
equivalence orbit of another (non-zero) random matrix $U A_1 V$ of
the same dimension,
since generically $A_0$ and $A_1$ clearly do not share the same singular
values. However, a random complex matrix $A_0$ is in fact typically
arbitrarily close to {\em a sum of two or more equivalence orbits of
independent random matrices}. This is shown in Fig.~\ref{fig:UAV_N2}
by a numerical example for $10\times 10$ complex square matrices,  where the
inset shows this does not hold for similarity orbits of random square
matrices. Interestingly, the findings hold independent of the dimensions and
explicitly include rectangular matrices as well as square matrices.

\bigskip
\noindent
{\bf Example 6}\\
For a single random complex square matrix $A_0 \in \IC^{10\times 10}$ we now
ask how close it typically is to the sum of $N=1,2,3,4,5,10$ equivalence
orbits $\sum_{j=1}^N U_j A_j V_j$, where the $A_j$ are independently chosen
random complex matrices $A_j \in \IC^{10\times 10}$. We compare the findings
with those of $N$ independent similarity orbits $\sum_{j=1}^N U_j A_j U^*_j$
and find the results of Fig.~\ref{fig:UAV_N2} underscoring Proposition 3.2.
%%%%%%%%%%%%%%%%%%%

\section{Unitary $t$"~Congruence}

In this section, we consider
$$\min \left\{\|\sum_{j=1}^N U_jA_j U_j^t - A_0\|: U_1, \dots, U_N \in U(n)
\hbox{ unitary}  \right\}$$
for given matrices $A_0, A_1, \dots, A_N$.
Sometimes, we can focus on special classes of matrices such as
symmetric matrices or skew-symmetric matrices. For symmetric
matrices or
skew-symmetric matrices,
the minimization problem
$$     \min\left\{ \|UAU^t - A_0\|: U \hbox{ unitary} \right\}$$
has an analytic solution; see \cite{MO}.
The problem is wide open even if $N = 2$.
Therefore, a computer algorithm will be most helpful in the theoretical development.
One may also consider whether we can have
$UAU^t + VBV^t = C$
for a symmetric $A$ and a skew-symmetric $B$.
In other words, we want to know whether one
can write $C$ as the sum of symmetric and skew-symmetric
matrices with prescribed singular values.
Of course, the problem for general matrices $A, B$ and $C$ is even
more challenging, and that is what we pursue by the numerical methods developed
in the next paragraph.

\subsection{Gradient Flows on Unitary $t$-Congruence Orbits}\label{sec:flow_congruence}
%%%%%%%%%%%%%%%
Again, the minimization task
\begin{equation}
\minover{U,V\in U(n)} ||UAU^t + VBV^t - C||_2^2\;,
\end{equation}
translates via
\begin{equation*}
||UAU^t + VBV^t - C||_2^2 \\
= ||A||_2^2 + ||B||_2^2 + ||C||_2^2 - 2 \Re \tr\left\{C^*(UAU^t + VBV^t) - UAU^t\; \Bar VB^* V^*\right\}
\end{equation*}
into maximising the function
\begin{equation}
F(U,V) := \Re f(U,V) := \Re \tr\left\{(UAU^t + VBV^t)\, C^* - UAU^t\; \Bar VB^*V^*\right\}\quad,
\end{equation}
where the differential reads (by virtue of the short-hand $\Tilde C:= C^* - \Bar VB^*V^*$)
\begin{eqnarray*}
D_U f(U)(\Omega U) &=& \tr\left\{((\Omega U) A U^t + UA(\Omega U)^t)(C^* - \Bar VB^*V^*)\right\} \\
           &=& \tr\left\{(\Omega U) A U^t \Tilde C\right\} + \tr\left\{ (UA(\Omega U)^t \Tilde C)^t\right\}\\
           &=& \tr\left\{( A U^t \Tilde C + A^tU^t\Tilde C^t) (\Omega U)\right\}\quad.
\end{eqnarray*}
%%%%%%%%%%%%
%%%%%%%%%%%%
From identifying
$
D_U f(U)\cdot (\Omega U) = \braket{\grad_U f(U)}{\Omega U} = \tr\left\{(\grad_U f(U))^* \Omega U\right\}
$
one finds
\begin{equation}
\grad_U f(U)  = (U A U^t \Tilde C + U A^tU^t\Tilde C^t)^* U
\end{equation}
so as to obtain for $F(U):= \Re f(U)$
\begin{equation}
\grad_U F(U) = -\big(UAU^t \Tilde C + U A^tU^t \Tilde C^t\big)_s\,U\quad.
\end{equation}
Again, taking the respective Riemannian exponentials
$\exp_U(\grad_U F(U))$ and $\exp_V(\grad_V F(V))$ thus gives
the slightly lengthy formula
\begin{equation}
U_{k+1} = \exp \left\{-\alpha_k \big( U_kAU_k^t (C^* - \Bar V_kB^*V_k^*)
                + U_kA^tU_k^t (C^* - \Bar V_kB^*V_k^*)^t \big)_s \right\}\;U_k%\\
%V_{k+1} &=& \exp \left\{-\beta_k \big( V_kBV_k^t (C^* - \Bar U_kA^*U_k^*)
%               + V_kB^tV_k^t (C^* - \Bar U_kA^*U_k^*)^t \big)_s \right\}\;V_k
\end{equation}
---and an analogous equation for $V_{k+1}$ by substituting $V$ for $U$ and $B$ for $A$---as
discretized solutions of the coupled gradient system
\begin{equation}
\dot U = \grad_U F(U,V) \quad\text{\rm and}\quad \dot V = \grad_V F(U,V) \;.
\end{equation}

\bigskip
Likewise, for higher sums of congruence orbits one finds
\begin{equation}
\min\left\{\| \sum_{j=1}^N U_j A_j U^t_j - A_0 \|:
U_1, \dots, U_N \in U(n)\; \hbox{ unitary} \right\}
\end{equation}
to be solved by the coupled system of flows ($j=1,2,\dots, N$)
\begin{equation}
U_{k+1}^{(j)} = \exp \left\{-\alpha_k^{(j)}
{\big( A_k^{(j)} A_{0jk}^* +  ( A_{0jk}^* A_k^{(j)})^t\big)}_s \right\}\;U_k^{(j)}\quad,
\end{equation}
where for short we set $A_k^{(j)} := U_k^{(j)}A_j{U_k^{^t(j)}}$ and
$A_{0jk}:= A_0 - \sum\limits_{\nu=1 \atop \nu\neq j}^N A_k^{(\nu)}$.

\section{Outlook: Non-Compact Groups}

For orbits $S(A)$ of matrices $A$ under the action of non-compact groups,
there are usually no good results for supremum or infinmum of the
quantity
$$\|X_0 - \sum_{j=1}^N X_j\|$$
with $X_j \in S(A_j)$ for $j = 0, 1, \dots, N$, for given
matrices $A_0, \dots, A_N$.

For example, for the invertible congruence orbit  of $A \in M_n$
$$S(A) = \left\{ S^*AS: S \in M_n \hbox{ is invertible} \right\},$$
we can let $S = rI$. Then
$$\|S^*A_0S - \sum_{j=1}^N S^*A_jS\|$$
converges to 0 or $\infty$ depending on
$r \rightarrow 0$ or $r \rightarrow \infty$.

Similarly, the same problems occur for
the equivalence orbit of $A \in M_n$
$$S(A) = \left\{SAT:  S, T \in M_n \hbox{ are invertible} \right\}.$$

For the similarity orbits, we have the following.

\begin{proposition} Suppose not all the matrices $A_0, \dots, A_N$ are scalar.
Then
$$\sup\|A_0 - \sum_{j=1}^N S_j^{-1}A_jS_j\| = \infty.$$
\end{proposition}

\it Proof. \rm
Suppose one of the matrices, say, $A_i$ is non-scalar.
Then there is $S_j$ such that $S_j^{-1}A_iS_j$ is in lower triangular form
with the $(2,1)$ entry equal to 1, and there are invertible matrices
$S_j$ such that $S_j^{-1}A_js_j$ is in upper triangular form
for other $j$. Let $D_r = \diag(r,1,1,\dots,1)$.
Then the sequence
$$(S_0D_r)^{-1}A_0 (S_0D_r)  - \sum_{j=1}^N(S_jD_r)^{-1}A_j (S_jD_r)$$
has unbounded $(2,1)$ entry as $r \rightarrow \infty$.
The conclusion follows. \qed

Determining
$$\inf\|A_0 - \sum_{j=1}^N S_j^{-1}A_jS_j\|$$
is more challenging.
Let us first consider two matrices $A, B \in M_n$.
We have the following.

\begin{proposition} Let $A, B\in M_n$. Then for any unitary similarity
invariant norm $\|\cdot\|$,
$$\|(\tr A-\tr B)I/n\| \le \|S^{-1}AS - T^{-1}BT\|$$
for any invertible $S$ and $T$.
\end{proposition}

\it Proof. \rm
Given two real vectors $x = (x_1, \dots, x_n),
y = (y_1, \dots, y_n)$, we say that
$x$ is weakly majorized by $y$, denoted by $x \prec_w y$ if the
sum of the $k$ largest entries of $x$ is not larger than that of
$y$ for $k = 1, \dots, n$. By the Ky Fan dominance theorem,
if $X = \diag(x_1, \dots, x_n)$ and $Y = \diag(y_1, \dots, y_n)$
are nonnegative matrices such that $(x_1, \dots, x_n) \prec_w (y_1, \dots,
y_n)$, then $\|X\| \le \|Y\|$ for any unitarily invariant norm $\|\cdot\|$.

Now, suppose $S^{-1}AS - T^{-1}BT$
has diagonal entries $d_1, \dots, d_n$ and singular values
$s_1, \dots, s_n$. Then
$$|\tr A - \tr B|=|\sum_{j=1}^n d_j|\le    \sum_{j=1}^n |d_j|.$$
Thus,
$$|\tr A - \tr B| (1, \dots, 1)/n \prec_w (|d_1|, \dots, |d_n|)
\prec_w (s_1, \dots, s_n).$$
It follows that
$$\|(\tr A - \tr B)I/n\|
\le \|\diag(|d_1|, \dots, |d_n|) \le \|\diag(s_1, \dots, s_n)\|
= \|S^{-1}AS - T^{-1}BT\|.$$
\vskip -.3in \qed

Can we always find invertible $S$ and $T$ such that
$$\|S^{-1}AS - T^{-1}BT\| = \|(\tr A-\tr B)I/n\|?$$
The answer is no, and we have the following.

\begin{proposition} Let $\|\cdot\|$ be a unitarily invariant norm
on $M_n$. Suppose $A \in M_n$ has eigenvalues
$a_1, \dots, a_n$, and $B = bI$. Then
$$\inf\left\{\|S^{-1}AS - B\|: S \in M_n \hbox{ is invertible} \right\}
= \|\diag(a_1-b, \dots, a_n-b)\|.$$
\end{proposition}

\it Proof. \rm
Suppose $S^{-1}AS - B$ has eigenvalues
$a_1 - b, \dots, a_n - b$, and singular values $s_1, \dots, s_n$.
Then the product of the $k$ largest entries
of the vector
$(|a_1 - b|, \dots, |a_n - b|)$ is not larger than
$(s_1, \dots, s_n)$ for $k = 1, \dots, n$.
It follows that
$$(|a_1 - b|, \dots, |a_n - b|) \prec_w
(s_1, \dots, s_n),$$
and hence
$$\|\diag(|a_1 - b|, \dots, |a_n - b|)\| \le
\|\diag(s_1, \dots, s_n)\| = \|S^{-1}AS-B\|.$$
Note that there is $S$ such that $S^{-1}(A-B)S$ is in upper triangular
Jordan form with diagonal entries
$a_1 - b, \dots, a_n - b$.  Let
$D_r = \diag(1, r, \dots, r^{n-1})$ for $r > 0$.
Then $(SD_r)^{-1}(A-B)(SD_r) \rightarrow \diag(a_1-b, \dots, a_n-b)$
and  $\|(SD_r)^{-1}(A-B)(SD_r)\| \rightarrow \|\diag(a_1-b, \dots, a_n-b)\|$
as $r \rightarrow 0$. So, we get the conclusion about the infinmum.
\qed

From the above result and proof, we see that if
$A$ has an eigenvalue $a$ with eigenspace of dimension $p$
and $B$ has an eigenvalue $b$ with eigenspace of dimension $q$
such that $p+q-n=r > 0$ then $S^{-1}AS-T^{-1}BT$ has an eigenvalue
$a-b$ of multiplicity at least $r$. The question is whether
we can write $A = aI_r \oplus A_1$ and
$B = bI_r \oplus B_1$ and show that
$$\inf \|S_1^{-1}A_1S_1 - T_1^{-1}B_1T_1\| = \|(\tr A_1-\tr B_1)I_{n-k}/(n-
k)\|.$$

\iffalse
Here is a simple case. If $A = \diag(1,1,-2)$ and $B \in M_3$ is a rank 2
nilpotent, can we show that $\inf\|S^{-1}AS - B \| = 0$?
\fi

\medskip\noindent
It is interesting to note that the following two quantities may be different.

\medskip
1) $\inf\left\{\|S^{-1}AS - T^{-1}BT\|: S \hbox{ is invertible} \right\}$.

\medskip
2) $\inf\left\{\|S^{-1}AS - B\|: S \hbox{ is invertible} \right\}$.

\medskip
For example, suppose $A = \diag(2,-1,-1)$ and 
$B=\(\begin{array}{ccc}0&1&0\\ 0&0&1\\ 0&0&0\end{array}\)$.  
Then there are invertible $S$ and $T$ such that
$$S^{-1}AS = \begin{pmatrix} 0&1&1\\ 1&0&1\\ 1&1&0\cr\end{pmatrix}
\quad \hbox{ and } \quad
T^{-1}BT = \begin{pmatrix} 0&1&1\\ 0&0&1\\ 0&0&0\cr\end{pmatrix}.$$
So, $C = S^{-1}AS-T^{-1}BT$ is a rank two nilpotent. 
Thus for any $\varepsilon > 
0$, there is an invertible $R_\varepsilon$ such that
$$R_\varepsilon^{-1}CR_\varepsilon =
\begin{pmatrix} 0&\varepsilon &0 \\ 0&0&\varepsilon\\ 0&0&0\cr\end{pmatrix}.$$
As a result,  
$$\|R_\varepsilon^{-1}S^{-1}ASR_\varepsilon 
- R_\varepsilon^{-1}T^{-1}BTR_\varepsilon\| \rightarrow 0 \quad
\hbox{ as } \quad \varepsilon \rightarrow 0.$$
So, the quantity in (1) equals zero.
On the other hand,
for every invertible $S$, we have
$$\|\(A-SBS^{-1}\)\(Se_1\)\|=\|A\(Se_1\)\|\ge \|Se_1\|$$
Therefore, $\inf\|A - SBS^{-1} \| \ge 1$.
So, we see that the quantities in (1) and (2) may be different.

\medskip\noindent
In connection to the above discussion, it is interesting to
study the following problem.

\medskip\noindent
1. Determine 
$$\inf\left\{\|S^{-1}AS-TBT^{-1}\|: S, T \hbox{ are invertible} \right\}$$
and characterize the matrix pairs $(A,B)$

\medskip\noindent
2. Determine 
$$\inf\left\{\|S^{-1}AS-B\|: S \hbox{ is invertible} \right\}$$
and characterize the matrix pairs $(A,B)$ 
attaining the infinmum if they exist.

\section{Conclusions}
We have treated the least-squares approximation problems by elements on the
{\em sum} of  various matrix orbits including unitary similarity, equivalence
and congruence.
Special attention has been paid to sums of unitary similarity orbits of a
Hermitian $A$ and a skew-Hermitian $B$, where theoretical
results have been obtained and shown to be consistent with
numerical findings. Further, new results on unitary equivalence
orbits have been obtained stimulated by numerical experiments.
are related to geometric arguments.

A general framework based on the gradient flows on matrix orbits
arising from Lie group actions has been developed to study the proposed
problems.
The gradient flows devised to this end extend the existing toolbox
(see e.g. \cite{Bloch94, Chu04})
by referring to sums of matrix orbits as summerized in Tab.~1.
This general approach can be used to treat many problems in theory
and applications. For instance,
flows on such sums of unitary similarity orbits
can also be envisaged as on unitaries taking a block-diagonal form,
and hence they relate to relative $C$~numerical ranges,
where the group action is restricted to a compact subgroup
$\mathbf K \subseteq SU(n)$ of the full unitary group \cite{SDHG07}.
Finally, first results on matrix orbits under non-compact group actions invite
further research.

%{\bf CK's Question:} Since the sup is $\infty$  and inf is 0 in some 
%cases and we know theorectically,  
%do we expect the gradient flow will lead to these sup and inf?
% internal answer TSH: YES we do; but we should not necessarily
% mention that in the paper too strongly.

%%%%%%%%%%%%%%%%%%%%%%%%%%%%%%%
%%%%%%%%%%%%%%%%%%%%%%%%%%%%%%%%%%%%%%%%%%%%%%%%%%%%%%%%%%%%%%%%%
\begin{table}[Ht!]
\caption{\label{tab:flows}
Summary of Least-Squares Approximations by
Matrix Orbits and Related Gradient Flows}
\begin{tabular}{lll}\\
\hline\\[-3mm]
type\; and\; objective && coupled gradient flows \\[1mm]
\hline\hline\\[-2.5mm]
{\em unitary similarity:}\\
{$\minover{U\in SU(n)} = \| \sum\limits_{j=1}^N U_j A_j U^*_j - A_0 \|$}
&&{$U_{k+1}^{(j)} = \exp\left\{-\alpha_k^{(j)} {[ A_k^{(j)}, A_{0jk}^*]}_s\right\}\;U_k^{(j)}$}\\[2mm]
&
&{\hspace{5mm} where \; $A_k^{(j)} := U_k^{(j)}A_j{U_k^{(j)}}^*$ %
	and $A_{0jk}:= A_0 - \sum\limits_{\nu=1 \atop \nu\neq j}^N A_k^{(\nu)}$} \\[7mm]
\hline\\[-2.5mm]
{\em unitary equivalence:}\\
{$ \minover{U,V \in SU(n)}  \|\sum\limits_{j=1}^N U_jA_jV_j - A_0\|$}
&&{$U_{k+1}^{(j)}=\exp\left\{-\alpha_k^{(j)} (U_k^{(j)} A_j V_k^{(j)} A_{0jk}^*)_s\right\}\;U_k^{(j)}$} \\[2mm]
&   &{$V_{k+1}^{(j)} = \exp\left\{-\beta_k^{(j)} (V_k^{(j)} A_{0jk}^* U_k^{(j)}A_j)_s\right\}\;V_k^{(j)}$} \\[2mm]
&   &{\hspace{5mm} where \;  $A_{0jk}:= A_0 - \sum\limits_{\nu=1\atop \nu\neq j}^N U_k^{(\nu)} A_\nu V_k^{(\nu)}$} \\[7mm]
\hline\\[-2.5mm]
{\em unitary congruence:}\\
{$ \minover{U\in SU(n)}  \| \sum\limits_{j=1}^N U_j A_j U^t_j - A_0 \|$}
&&{$U_{k+1}^{(j)} = \exp\left\{-\alpha_k^{(j)}
{\big( A_k^{(j)} A_{0jk}^* +  ( A_{0jk}^* A_k^{(j)})^t\big)}_s\right\}\;U_k^{(j)}$}\\[2mm]
&   &{\hspace{5mm} where \; $A_k^{(j)} := U_k^{(j)}A_j{U_k^{(j)}}^t$ %
	and $A_{0jk}:= A_0 - \sum\limits_{\nu=1 \atop \nu\neq j}^N A_k^{(\nu)}$} \\[7mm]
%\hline\\[-2.5mm]
%{$ \minover{U\in SU(n)} \| \sum\limits_{j=1}^N U_j A_j \bar U_j - A_0 \|$}
%&&{$U_{k+1}^{(j)} = \exp\left\{-\alpha_k^{(j)}
%{\big( A_k^{(j)} A_{0jk}^* - U_k^{(j)} ( A_{0jk}^* A_k^{(j)})^t \;U_k^{*(j)}\big)}_s\right\}\;U_k^{(j)}$}\\[2mm]
%&   &{\hspace{5mm} where \; $A_k^{(j)} := U_k^{(j)}A_j{\bar U_k^{(j)}}$ %
%	and $A_{0jk}:= A_0 - \sum\limits_{\nu=1 \atop \nu\neq j}^N A_k^{(\nu)}$} \\[2mm]
\hline\\[-2.5mm]
%%%%%%
%%
\end{tabular}
\end{table}
%%%%%%%%%%%%%%%%%%%%%%%%%%%%%%%%%%%%%%%%%%%%%%%%%%%%%%%%%%%%

%%%%%%%%%%%%%%%%%%%%%%%%%%%%%%%
%%%%%%%%%%%%%%%%%%%%%%%%%%%%%%%

\section{Further Research}

In order to avoid the search in our algorithms is terminated in local extrema,
one has to ensure to choose a sufficiently large
set of random unitaries distributed according to the Haar measure.
Actually, one knows there are commutation properties
at the critical points. It would be nice to find a more
efficient method to choose starting points for the search,
and prove theorems ensuring that the absolute minimum will
be reached from one of these starting points using our
algorithms.

Our discussion focused on orbits of matrices under
actions of compact groups. We can consider other orbits
under actions of non-compact groups. Here are some examples
for $S,T \in SL(n,\IC)$:

(e) the general similarity orbit of a square matrix $A$ is the set
of matrices of the form $SAS^{-1}$,

(f) the equivalence orbit of a rectangular matrix $A$ is the set
of matrices of the form $SAT$,

(g) the $*$-congruence orbit of a complex square matrix $A$ is the set of
matrices of the form $SAS^*$,

(h) the $t$-congruence orbit of a square matrix $A$ is the set of matrices
of the form $SAS^t$.

\noindent
However, the fact that $GL(n,\IC)$ and $SL(n,\IC)$ are just {\em locally compact}
entails there is no Haar measure and consequently no {\em bi-invariant} metric
on the tangent spaces, but only left {\em or} right-invariant metrics. Hence
the Hilbert-Schmidt scalar product $\braket{B}{A}=\tr\left\{B^* A\right\}$ has to be treated
with care, in particular since we are interested in the complex domain.
Moreover, while in compact Lie groups the exponential map is surjective and
geodesic \cite{Arv03}, in {\em locally compact} Lie groups, it is generically
neither surjective nor geodesic. It is for these reasons that
devising gradient flows in locally compact Lie groups is the subject of a follow-up study.

%\medskip
%There are other related question such as studying
%
%(i) determine $\min \|e^{iH}e^{iG} - e^{i(U^*HU+V^*GV)}\|$.\\
%
%Related to the $UAU^* + VBV^* = C$ problem, we consider
%$$\min \left\{\|e^{iUAU^*}e^{iVBV^*} - e^{iC}\|: U, V \hbox{ are unitary}
% \right\}$$
%for given Hermitian $A, B, C$.\\
%{\em Chi-Kwong: for the moment, we have to skip the paragraph (i)
% due to bad numerical conditioning that cannot be eliminated
% in a stable manner by conventional scaling and squaring.}

\end{document}